\documentclass{article}

\usepackage{graphicx,amsmath}

\begin{document}

\title{Multi-criteria optimization methods in radiation therapy planning: a review of technologies and directions}

\author{David Craft \\ ~ \\ Massachusetts General Hospital, Boston, MA \\ ~ \\}

\maketitle

\abstract{We review the field of multi-criteria optimization for radiation therapy treatment planning. Special attention is given to the technique known as Pareto surface navigation, which allows physicians and treatment planners to interactively navigate through treatment planning options to get an understanding of the tradeoffs (dose to the target versus over-dosing of important nearby organs) involved in each patient's plan. We also describe goal programming and prioritized optimization, two other methods designed to handle multiple conflicting objectives.  Issues related to nonconvexities, both in terms of dosimetric goals and the fact that the mapping from controllable hardware parameters to patient doses is usually nonconvex, are discussed at length since nonconvexities have a large impact on practical solution techniques for Pareto surface construction and navigation. 
A general planning strategy is recommended which handles the issue of nonconvexity by first finding an ideal Pareto surface with radiation delivered from many preset angles. This can be cast as a convex optimization problem.  Once a high quality solution is selected from the Pareto surface, a sparse version (which can mean fewer beams, fewer segments, less leaf travel for arc therapy techniques, etc.) is obtained using an appropriate sparsification heuristic. We end by discussing issues of efficiency regarding the planning and the delivery of radiation therapy.}

\section{Overview of radiation therapy}
Surgery, radiation, and chemotherapy are the three most common types of cancer treatment. Radiation treatment, used in approximately two-thirds of all cancer treatments, is a technically challenging modality involving ever-improving medical imaging and precise radiation production and beam shaping. This review focuses on multi-criteria optimization techniques used in the design of radiation therapy plans.

Radiation is considered a localized therapy in comparison to chemotherapy, which is systemic. However, radiation dose cannot be delivered to only the target: there is inevitably a dose fall off that can and does irradiate nearby tissues and organs, which gives rise to the main planning challenge in radiation treatment: balancing high levels of radiation dose to targets with detrimental effects of dose to nearby non-cancerous tissues.

\begin{figure}[h!t]
\label{hierarchy}
\centering
\includegraphics[trim=0 200 0 100,clip,width=\textwidth]{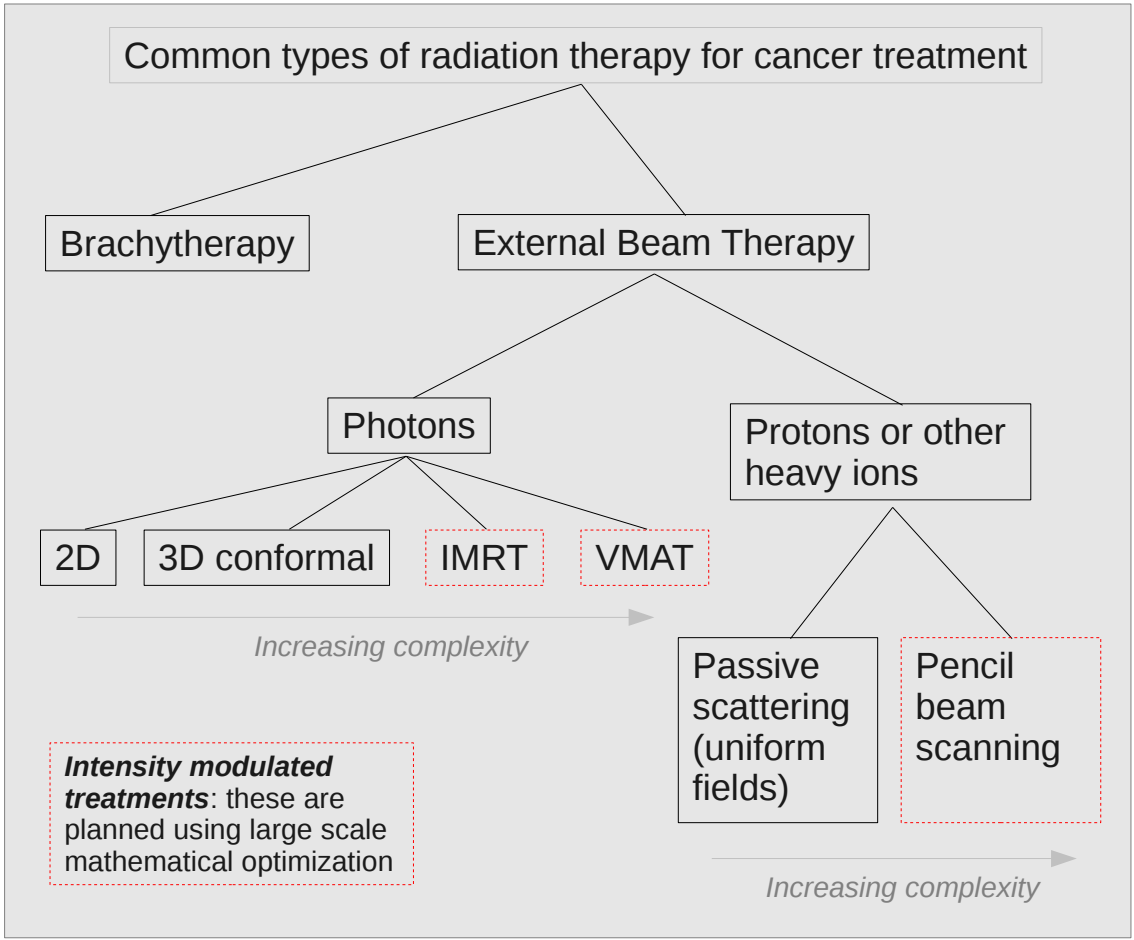}
\caption{Hierarchy of the most common radiation delivery types for cancer treatment. The highlighted boxes are the types of treatment that we focus on in this review.}
\end{figure}

Radiation therapy comes in two main forms: internal (brachytherapy) and external. External beam therapy is more common. Brachytherapy, which is highly localized, consists of implanting radioactive seeds into the tumor and therefore requires a surgical procedure. External beam therapy does not require a radioactive seed source, is applicable to a wider range of disease sites, and is considered equally effective. Here we will focus on planning external beam radiation therapy although many of the principles, in particular the notion of dealing with the tradeoff between effectively irradiating the tumor and not over irradiating the nearby healthy tissues, apply equally well to brachytherapy. We look at forms of external beam radiation therapy which require large scale computerized optimization algorithms (see Figure \ref{hierarchy} for a hierarchical sketch of the most common types of radiation therapy) and consider the set of approaches which fall under the umbrella of multi-criteria optimization, a general class of techniques designed to handle conflicting goals of a design problem.

\subsection{Intensity modulated treatments}

Intensity modulation in external beam radiation therapy refers to the use of treatment fields that are non-uniform, as opposed to flat fields. Intensity modulation techniques entered the field of radiation therapy in the mid 1990s when computer algorithms and control systems were made available which could control a multi-leaf collimator (MLC) placed in the beam path to block portions of the photon beam for differing amounts of time, thus creating arbitrarily shaped beam profiles \cite{Webb2003}. The additional degrees of freedom of non-uniform beams gives the chance to tailor the dose distribution in the patient to a much greater extent, allowing the creation of treatment plans that carve dose away from critical structures and/or give additional dose to particularly resistant regions. Controlling the large number of degrees of freedom however demands the use of optimization algorithms since exploring the space of potential plans manually would be far too impractical.

Pencil beam scanning proton therapy is emerging as another form of intensity modulated radiation therapy in addition to more common X-ray based treatments. Proton beams have the advantage of depositing the peak dose far beyond the patient's skin surface, which is obviously beneficial to tumors occuring deep inside the body. Mathematically, photons (high energy X-rays) and protons give rise to the same type of optimization problem and therefore we will not distinguish between the two modalities herein. 

Our basic optimization model is given next. We assume the patient computed tomography (CT) scan has been divided into a grid of voxels which are typically cubes with side length between 2 and 4 mm. Dose to each voxel are the primary quantities of interest, which we write as a dose vector $d$, where the length is equal to the total number of patient voxels, order of $10^6$. The control variables are the beamlet intensities (which result from dividing up each radiation field into a regular grid with elements typically of size 5 mm $\times$ 5 mm or 1 cm $\times$ 1 cm), written as a vector $x$, length on the order of $10^4$. We assume the mapping from beamlet intensities to voxel doses is given by the matrix $D$. $D$ is sometimes refered to as the dose-influence matrix. The basic radiation therapy optimization problem is:
\begin{eqnarray}
\hbox{ minimize~} & f(d) \nonumber \\
~&  Dx = d \nonumber \\
~&  d \in C \nonumber \\
~&  x \ge 0,
\label{basic}
\end{eqnarray}
\noindent If $f$ is a convex function and $C$ is a convex set, this is a convex optimization problem. In general however, particularly in commercial treatment planning systems which use dose-volume constraints and objectives (i.e. control the number of voxels above or below a certain dose level), $f$ and $C$ are non-convex, which makes the optimization theoretically much harder, since a local minima is not guaranteed to be a global minima. In practice, local minima are typically considered good solutions, and there are theoretical reasons why they might in fact be close to global optimal solutions, see Section \ref{nonconvexsection} for more details. Therefore dose-volume constraints are not thought to be a problem in practical treatment planning \cite{localMinima}. 

While the optimization problem above is typically large, it is not the computation time that is the main challenge in intensity modulated treatment design. Instead, it is dealing with the natural target dosing versus normal tissue sparing tradeoff. With intensity modulation, the variety of possible dose distributions for a given patient has skyrocketed, which makes the search for the ``best plan'' time-consuming. The typical process for a difficult case in the clinic is an iterative one between the physician and the treatment planner where the physician states a set of desires, the planner attempts to come up with a treatment plan that matches those desires--by tweaking the $f$ function and the constraints $C$--and they iterate until either a satisfactory plan is created or until there is no time left for further plan alterations. This process is inefficient and often leads to suboptimal plans. The use of sensitivity analysis \cite{sens2,sens3}--analyzing how small changes in constraints change the objective value--can be useful but with a large number of constraints and objectives that often do not have a clear clinical meaning (e.g. sum of the quadratic deviations from the prescription dose) and the fact that sensitivity information is by definition local information (i.e. only applies to small changes in the constraint values) it is of limited use in clinical practice. Indeed, most IMRT treatment planning systems do not offer sensitivity analysis.

\subsection{What makes a good treatment plan?}

Since the dose distribution to an organ is typically not uniform over that organ, quantitatively dealing with the dose tradeoffs must begin by succinctly describing the dose to an organ or a target. Functions which are currently used include the mean dose, the maximum dose, or the fraction of volume of a structure that receives more or less than a certain specified dose level. What would be more useful would be functions that mapped dose distributions into toxicity or cure outcomes for that structure. Such models, termed NTCP (normal tissue complication probability) and TCP (tumor control probability) are often discussed in the community but rarely used to design clinical plans, due to the lack of consolidated clinical data to estimate these functions \cite{nahum}. This complicates the design of good plans since without such descriptors, a fair amount of hand waving enters the picture in terms of discussing the relevant merits of one plan over the next. Indeed, even with reliable NTCP and TCP models, it is a subjective matter--and dependent on patient preferences--to say whether a plan with high TCP and NTCPs is more preferable to a plan with lower TCP and NTCPs. We also note that there are multiple organs to consider in any single treatment (cancers of the head and neck region can easily have over 10), and each organ might have more than one potential complication type, and hence NTCP function. Therefore, the tradeoff is not two-dimensional: TCP vs NTCP, but multi-dimensional: TCP vs NTCP$_1$ vs NTCP$_2$ vs NTCP$_3$, etc. Finally, the idea of reducing a dose distribution to a set of individual TCP and NTCP functions is itself a simplification because this implicitly assumes organ independence. In reality, the dose to one organ might influence the NTCP to a different organ. For example it has been shown that the dose delivered to the heart is a predictive factor for radiation pneumonitis of the lung \cite{lungheart}, but such interdependencies are rarely taken into account in treatment plan design. 

This short discussion provides some backdrop to the use of multi-criteria optimization (MCO) for treatment planning. Clearly a treatment plan cannot be described by a scalar value: a treatment plan good at say sparing the liver, might give a too aggressive dose to the stomach. Managing these inherent tradeoffs, and presenting them clearly to an informed decision maker, is the goal of multi-criteria optimization. At this point, we will assume we have in our hands mathematical functions such as TCP and NTCPs for the various involved organs (clinically, mean dose, maximum dose, dose-volume criteria, etc. are used, but it is clearer for exposition purposes to assume we have TCP and NTCP functions), and proceed with the main topic of this review: managing the complex high dimensional tradeoff amongst these entities in radiation therapy planning. In the discussion section we will return to the important reality that such functions are not known, and are patient/genome/health status dependent. 

\section{Three types of MCO: goal programming, prioritized optimization, and Pareto surface navigation}

It is possible to frame the radiation therapy planning problem as a single objective optimization problem: for example, maximize TCP subject to the constraints that the NTCP functions for the relevant organs are not too large, for example $\le$ 3\%. Such formulations can be useful and can lead to good treatment plans. However, after viewing the result from such an optimization, the physician or planner may be curious to see what type of improvement in TCP would result from raising the NTCP tolerances by a couple of percent. If there are multiple structures to consider, such what-if scenarios make for very time-consuming treatment planning when not handled well.

{\bf Goal programming} casts the treatment planning problem as a set of goals. For example, perhaps the doctor specifies that the TCP $\ge$ 90\% and all NTCP functions $\le$ 5\%. The goal programming method tries to find a plan which satisfies a set of such constraints. Two things can happen at this point: either a plan can be found which satisfies these goals, or not. In the first case, one could naturally ask: can one do better, and in the second case one would ask, what goals should I sacrifice in order to achieve an acceptable plan. These two scenarios indicate a deficiency of the goal programming approach for radiation therapy. However, given the clinical reality that physicians often specify their desires in such language, the goal programming framework is relevant to modern treatment planning.

For concreteness, assume we are trying to minimize given functions $f_j$. For example, these could be dose-volume functions to reduce the number of voxels which violate a dose level, or $f_j$ could be the NTCP for structure $j$. Let $g_j$ denote the goal that we are trying to reach. The goal programming radiation treatment planning problem can then be written as:
\begin{eqnarray}
\hbox{ minimize~} & \sum_j s_j \nonumber \\
~&  Dx = d \nonumber \\
~&  s_j \ge f_j(d) - g_j \nonumber \\
~&  s_j \ge 0 \nonumber \\
~&  d \in C, x \ge 0
\label{gp}
\end{eqnarray}
\noindent If $f_j(d) \ge g_j$, that is, the goal has not been met, $s_j$ will take on the value of the difference, or `slack'. If the goal has been met, $s_j$ will take on the value 0, hence no penalty in the objective functon is incurred. If all goals are met the value of the objective function will be 0, and there is no incentive built into this optimization formulation to lower the functions $f_j$ any further. 

Some physicians, rather than giving a general set of goals, will prioritize the goals. This prescription style leads naturally to {\bf prioritized optimization}, also called lexicographic optimization, and it is implemented by running a sequence of optimizations \cite{wilkens, falkinger, jee, clark}. It is in general useful to distinguish between objectives and constraints, something which is often not done in the radiation treatment planning community. In prioritized optimization it is critical to do so. Constraints are non-negotiable: they are the highest priority. Objectives, as used in the mathematical optimization community, have a sense (either maximize or minimize) but do not have a hard level they are supposed to reach. So, a prioritized optimization approach to designing a treatment plan with a maximum spinal cord dose of 45 Gy as a (non-negotiable) constraint might look as follows. Step 1: maximize TCP subject to spinal cord dose $\le$ 45 Gy. Let us assume that the TCP obtained from this optimization is 98\%. This value is turned into constraints in subsequent steps, including a relaxation factor, also known as a slip factor. Step 2 could then be: minimize mean lung dose subject to spinal cord dose $\le 45$ Gy and TCP $\ge 97$\%. Here an absolute slip factor of 1\% has been used. Lower priority goals would be handled in subsequent steps, each one of which will gather the previously obtained objective values and use them, with slip factors, as constraints. Formally, prioritized optimization proceeds as follows. Let $d \in C$ denote the non-negotiable constraints, and let $f_1$ be the first priority objective, $f_2$ be the second, etc. Step 1:

\begin{eqnarray}
\hbox{ minimize~} & f_1(d) \nonumber \\
~&  Dx = d,~~d \in C,~~  x \ge 0
\label{po1}
\end{eqnarray}
\noindent Then, letting $f_1^*$ denote the optimal objective value for the first step, the next optimzation proceeds as:
\begin{eqnarray}
\hbox{ minimize~} & f_2(d) \nonumber \\
~&  f_1(d) \le f_1^*+\epsilon  \nonumber \\
~&  Dx = d,~~d \in C,~~  x \ge 0
\label{po2}
\end{eqnarray}
Assuming a total of 3 priority levels, the optimization which yields the final treatment plan would be:
\begin{eqnarray}
\hbox{ minimize~} & f_3(d) \nonumber \\
~&  f_1(d) \le f_1^*+\epsilon  \nonumber \\
~&  f_2(d) \le f_2^*+\epsilon  \nonumber \\
~&  Dx = d,~~d \in C,~~  x \ge 0
\label{po3}
\end{eqnarray}

Both goal programming and prioritized optimization provide only a single plan at the end of the process, without a natural way to explore the plan tradeoffs. {\bf Pareto surface navigation} is an MCO technique which puts the interactive exploration of the dosimetric tradeoffs at the front and center. For a given set of objectives and constraints, a treatment plan is Pareto optimal if it satisfies all of the constraints and none of the objective values can be improved upon without worsening at least one of the other objectives. These are the plans that should be considered for patient treatment. The set of Pareto optimal treatment plans (termed the Pareto surface) for continuous optimization problems is typically an infinite set. Finding efficient ways to approximate it and intuitive ways to explore it are important research and development areas. The Pareto optimization problem is written as:
\begin{eqnarray}
\hbox{ minimize~} & \{ f_1(d), f_2(d), \ldots f_N(d) \}  \nonumber \\
~&  Dx = d \nonumber \\
~&  d \in C \nonumber \\
~&  x \ge 0,
\label{pareto}
\end{eqnarray}
\noindent where $N$ is the number of objective functions.
  
The goal of the Pareto surface navigation approach is to put the planning process in the hands of the decision maker most qualified to make the tradeoff decision, which in radiation therapy is the physician. Not all research and clinical groups agree on the role, if any, that Pareto surface navigation should play in a clinical setting. One argument against it is that treatment planning standards should be in place, thus promoting quality and reproducibility. In this case, either a goal programming approach with fixed objective function weights, or a defined prioritized optimization procedure, would be more suitable. Proponents of the Pareto surface navigation approach will argue that since each patient comes with distinct a cancer shape and nearby organ geometries, and a unique health-history profile, a manual exploration of the dosimetric tradeoffs performed by the treating physician is a sensible approach, if it can be done in a time-effective manner. Given the variety of cancer types that a modern radiation therapy clinic sees (including protocol-driven clinical trials, standard disease sites such as prostate and breast, and patient specific difficult targets such as anal cancers, head and neck, and whole abdominal radiation), it seems perfectly reasonable to have each of the three MCO approaches described above as part of a radiation therapy software suite, where clinics can decide which technique makes the most sense for each site. 

For the remainder of the technical portion of this review, we will focus our attention on Pareto surface navigation since that is the most difficult from a mathematical perspective, and also the most general of the techniques.

\section{Computational aspects of Pareto surface navigation}

\subparagraph{Computing the Pareto surface}

The first step of computing a Pareto surface is to define the objectives and the constraints. This step is somewhat neglected in the field since it does not pose a challenging mathematical problem, but the choice of formulation has dramatic downstream effects and should not be ignored by researchers or practioners.

Constraints are chosen to limit the range of the Pareto surface. If constraints are chosen too aggressively, the result might be that there are no feasible plans, and thus no Pareto surface. On the other hand, if constraints are chosen too loosely, the extent of the Pareto surface might be too large, thus requiring too many plans to approximate it. Ideally the Pareto surface contains a variety of plans that are all potential treatment candidates. Constraints chosen should be non-negotiable (this is in fact the definition of a constraint). Negotiable requests should be framed as objectives. The number of objectives chosen should be such that all relevant tradeoffs are exposed, but no more, since computational burden rises with the number of objectives, although perhaps not as fast as one might anticipate from a purely geometrical argument \cite{craft-hmp, spalke}. 

There are two main methods for computing Pareto optimal plans. For a comprehensive technical review as of 2005, see \cite{ruzika}.  The constraint method, also called the epsilon constraint method, constrains all but one of the objectives to achievable levels, and then minimizes the remaining objective. Variations on this method use objective function constraining in other ways but all operate on the same fundamental idea: to use a single minimization with some additional constraints to find a point in a particular region of the Pareto surface, e.g. \cite{lo,messac,wei-proj}, and repeat with different constraint values to find a set of points on different regions of the Pareto surface. Constraint methods are beneficial because methods in that class exist that can generate points anywhere on the Pareto surface even when the formulation is non-convex. Weighted sum methods on the other hand are attractive for convex formulations but are theoretically deficient in other settings. These methods minimize a non-negative weighted sum of the objectives, and are attractive because they work well with sandwich algorithms, which produce inner and outer approximations to the Pareto surface, allowing the approximation error to be controlled \cite{rennen, bokrantz}.

\subparagraph{Navigating the Pareto surface}

By navigation, we mean real time exploration of the continuous approximation of the Pareto surface. We start with the assumption that the Pareto surface has been approximated by a finite set of representative points (i.e. treatment plans). We also make the assumption, fully in line with the above convexity assumptions, that convex combinations (non-negative weighted sums of the plans where the weights add to 1) are valid plans. 

Compared to computing the Pareto surface, far less work has been published on the real-time navigation (exploration) of Pareto surfaces.  For 2D and 3D surfaces, one can visualize the surfaces fully and thereby graphically explore them. The vast majority of papers on computing Pareto surfaces are methods for 2D or 3D but which fail or are not intended for higher dimensions, which is one reason that so little work exists on navigating higher dimensional surfaces. A notable exception is \cite{monz}, which relies on linear programming to compute convex combinations of the pre-computed Pareto optimal plans which correspond to a user's requests as made through slider bars and adjustable objective function bounds. Other high dimensional navigation strategies that use similar linear programming methods include the commercially available multi-criteria IMRT treatment planning software RayStation (RaySearch Laboratories, Stockholm, Sweden) and Astroid, the in-house proton planning system developed at the Massachusetts General Hospital \cite{astroid, craftmonz}. A method based on viewing 2D tradeoff curves derived from projecting the higher dimensional data onto two selected objective axes is described in \cite{delivnav}.

All of these navigation methods depend on convex combinations of pre-computed plans, and therefore are designed with fluence map based plans in mind. For proton beam scanning this is a very good representation of reality, but for photon based delivery with an MLC fluence based plans need to be converted into MLC deliverable plans. Navigating on a Pareto surface of already deliverable plans is problematic because combining several MLC segmented plans results in a plan that uses all of the MLC segments of the underlying plans,  and too many segments make a plan clinically unattractive because of the unnecessarily long delivery times. The underlying problem here is the nonconvexity of cardinality constraints (in this case, the allowed number of MLC segments). In the next section we address this and other nonconvexity issues in radiation therapy planning.


\section{Nonconvexities}
\label{nonconvexsection}
Two types of nonconvexities occur in radiation therapy optimization: dosimetric and hardware derived. The primary example of dosimetric nonconvexities are dose-volume constraints. A constraint is nonconvex when a weighted average of two solutions which each satisfy the constraint may fail to satisfy the constraint. Dose-volume constraints have come to dominate the field of radiation therapy for historical reasons. While technically they are nonconvex, other planning constraints and the physical properties of radiation beams can make the nonconvexities disappear or be negligible. For example, consider the following dose-volume constraint: volume of stomach receiving 40 Gy or more should be less than 10\%. For say pancreatic cancer, most viable plans which dose the pancreas to 50 Gy (the prescription dose) will achieve this dose-volume constraint in a similar way: carve out the dose from stomach leaving only a small amount of volume of the stomach, the volume closest to the pancreas, with doses above 40 Gy. It is highly likely that the dosimetrically averaged plan created from two plans that each fulfill this constraint will still fulfill the constraint, due to the physical restrictions on how such plans can look.  Therefore, while theoretically leading to nonconvexity, dose-volume constraints are not necessarily nonconvex. Furthermore, it is likely that if two plans are each acceptable by a physician, then so would a dosimetrically averaged plan. For these reasons, we will not consider the potential nonconvexities that arise from dose-volume constraints.

Hardware derived nonconvexities arise from the practical desire to use a small number of beam angles and/or IMRT multi-leaf collimator segments when creating treatment plans. Consider the following modification of the basic radiation therapy problem:
\begin{eqnarray}
\hbox{ minimize~} & f(d) \nonumber \\
~&  Dx = d \nonumber \\
~&  d \in C \nonumber \\
~&  x \ge 0,~~~x \in S
\label{nonconvex}
\end{eqnarray}
\noindent Here, $S$ represents a (likely nonconvex) set of vaid delivery options. For example, $S$ could be the set of all five beam plans, or all plans formed of 100 or fewer multi-leaf collimator segments. These are truly nonconvex restrictions because the ``average'' of two plans, each using five beams, is a plan with up to ten beams. These types of cardinality restrictions necessitate global search techniques to find the optimal solutions. Typically strategies for these problems which would result in provable optima, such as mixed integer programming, are not computationally feasible without large computer infrastructure.

We will discuss in more detail three particular instances of hardware derived nonconvexities: beam angle optimization (BAO), direct aperture optimization (DAO), and volumetric modulated arc therapy (VMAT), in general and in regards to Pareto surface navigation.

The challenge for BAO regarding Pareto surface navigation is that the optimal beam angle set will vary across the Pareto surface. Although attempts have been made to characterize the BAO problem as a continuous problem in beam angle space \cite{craft-bao}, most BAO algorithms assume a discrete sampling of the beam angle space, thus a finite computation of the dose-influence matrix $D$. In this setting the Pareto surface will be a patched surface, different patches corresponding to different beam angle configurations for that subregion of the Pareto space.  The navigation of such a structure has been discussed \cite{craftmonz}, but the computation of it has not been addressed. This is partly due to the fact that even a single criteria BAO problem is very difficult to solve to proveable optimality. If one considers a 10 degree grid, assuming coplanar beams spaced evenly around the patient, and a seven beam plan, the number of possible beam angle sets to consider is ${36 \choose 7} \approx 8.3$ million. This makes a direct search prohibitively expensive, although mixed integer optimization techniques \cite{evasMIPBAO} and other global search strategies \cite{aleman2009, BAOmdandersonWANG} have been tried with some success, but require large computational infrastructure.  Due to the large search space, BAO methods are most often heuristics, see e.g. \cite{icycle, valentina} for approaches and reviews of the literature. In commercial treatment planning systems, beam selection is still done manually.

DAO, similar to BAO, attempts to find a sparse optimal solution to the IMRT planning problem \cite{dao,dao2}.  In the case of BAO, the sparsity is regarding the number of beam angles used. In DAO, the sparsity is regarding the number of MLC segment shapes used to deliver the radiation. Thus, DAO applies to the IMRT delivery technique known as step-and-shoot, and does not apply to sliding window (dynamic) delivery. The benefit of using a DAO approach for IMRT treatment planning is that one avoids the lossy conversion step from idealized fluence maps into MLC segment shapes when using the more traditional two-step IMRT optimization method. The downside is that the DAO problem, because it is a cardinality restriction problem, is nonconvex. Similar to the BAO problem, the DAO problem in the Pareto navigation  world has the same feature that different regions of the Pareto surface will have different optimal MLC segments. Two complimentary approaches for dealing with this are discussed in \cite{delivnav, daomco}. In \cite{daomco}, the researchers construct a 2D Pareto surface out of a single set of apertures by sequentially adding apertures which improve the entire surface thus far computed or a specific region. In \cite{delivnav} the researchers show that for high dimensional Pareto surfaces, the number of plans needed for the convex combinations formed during navigation can be restricted, thus promoting the idea that different parts of the Pareto surface can have different active aperture sets and that combining plans during navigation may not necessarily make the total number of apertures used for navigated plans too large. More sophisticated methods where apertures enter and leave the active set and morph during navigation have not appeared in the literature yet.

An interesting new intensity modulation modality that offers the chance to avoid the difficulties of both beam angle selection and step and shoot aperture selection is VMAT. VMAT delivers radiation continuously as the gantry rotates around the patient and the MLC leaves move back and forth across the field. Ideally, the gantry may slow down where extra modulation and/or dose is needed, but constant gantry speed VMAT is another choice. The nonconvexity of VMAT arises because the dose to a voxel is a nonconvex function of the underlying control variables, leaf and gantry positions as a function of time. Nevertheless, VMAT is particularly amenable to convex relaxations, for example the relaxation of full intensity modulation (where the control variables become beamlet weights and hence the problem is convex again) at a fine gantry spacing, say every 5 or 10 degrees \cite{vmerge, networkvmerge, rasmus-mco-vmat}. In those works, the following Pareto planning strategy emerged: 1) The user navigates a Pareto surface derived from a many-beam (convex) IMRT formulation. 2) After selecting a plan with the desired dosimetric tradeoffs, a true VMAT plan is produced using a heuristic which converts the fluence based plan to a deliverable VMAT plan, while deviating at most a user defined $\epsilon$ amount away from the ideal plan. Using this as a general purpose radiotherapy planning approach is discussed next.

\section{Deriving sparse near optimal plans from ideal (many-beam) plans} 

Sparsifying (i.e. simplifying) a dense solution to make it efficiently deliverable works well in many radiotherapy situations, but not all. It works well when the sparse solution sought is not so sparse that optimality is badly affected. The $l_0$-norm of a vector $x$, also known as the cardinality of $x$, is equal to the number of non-zero components of $x$. So, if $x$ is a vector whose components are the weights of all possible step-and-shoot apertures, then finding a solution with minimum number of apertures is a min $l_0$-norm of $x$ problem, which is in general a nonconvex, NP hard problem \cite{mangasarian}. The $l_1$-norm has been popularized as a convex substitute for the $l_0$ norm minimization problems in the field known as compressed sensing \cite{compress}, but the theory of compressed sensing applies to a situation where there is an underdetermined linear system and highly sparse solutions exist. The radiation therapy problem does not take the form of an underdetermined linear system, so compressive sensing ideas do not apply. Nevertheless, one may ask if a similar sparisifaction theory might exist for radiation therapy problems.  

The answer is likely muddier than one would hope: up to a certain point, plans can be vastly simplified. In IMRT, for example, there is no need for 100 beams, although at some point incremental plan degradation occurs with each beam or segment taken away, and this is a true tradeoff that requires expensive optimization to characterize. Initial simplification is possible \cite{zhu2009, dassimrt, craft-spg}, but if one wants to find the optimal seven beam IMRT plan for a complex head and neck case, this is truly a nonconvex problem with lots of distant local minima, and sparsification heuristics related to $l_0$ minimization \cite{candes, mangasarian}, if pushed to find such a sparse solution, will be pushed towards one of the local minima, but not necessarily the best one. In these cases, to find the true optimal solution, mixed integer programming or other global search strategies must be employed \cite{evasMIPBAO, BAOmdandersonWANG}.  To find good but not proveably optimal solutions to the small number of beams or apertures problems, many different heuristic approaches have been shown to work, including adding only beams that score favorably regarding metrics based on the beam's-eye-view of target coverage versus organ sparing \cite{pugachev02}, or metrics based on gradient information containing quantification of objective function improvement \cite{icycle}. Starting with a large collection of beams or aperatures and then successively removing them if they contribute little to the solution is also a promising approach. 


The vast simplification that is usually possible--that is, cutting down the number of beams, or segments, or the number of times the leaves travel across the collimator head in sliding window VMAT delivery--makes the following general approach, the same approach as discussed for VMAT optimization, appealing: 

\bigskip

\fbox{
  \parbox{13.4cm}{
     {\bf Ideal Pareto Surface Technique} \\ Begin with enough beams, segments, etc., and create a so-called \emph{ideal Pareto surface}. Navigate this surface to determine an \emph{ideal achievable dose distribution}. Lastly, reproduce the selected plan to a user selected $\epsilon$ level of closeness with a sparsification technique (some of which were discussed above) that is designed for the particular modality of interest.
  }
}

\bigskip

This approach is logical for radiation therapy planning problems in that its first focus is on achieving a plan of good balance between target coverage and healthy organ sparing. Only after this high quality plan is selected does it address the issue of plan deliverability. This works in the context of radiation therapy because there we know that much sparsification can be done with minimal impact to the dose distribution quality. Furthermore, delivery improvements that are here or on the near horizon, including consolidated beam delivery (see Section \ref{conc}), faster leaf speeds, and faster energy layer switching in proton therapy beam scanning, will make even solutions that by today's standards would not be considered sparse fast enough for clinical delivery.

\section{Conclusions: efficiency in treatment planning and delivery, and the role of MCO}
\label{conc}





Multi-criteria optimization has the potential to greatly improve treatment planning processes and treatment plans in radiation oncology, but there are several hurdles that must be overcome to reap the full benefits. As the name implies, multi-criteria optimization implicitly assumes that you have multiple criteria defined through which one can judge the quality of a treatment plan, but there are no widely agreed upon metrics for these. Moving to models that are directly tied to probabilities of certain outcomes--termed biological models in the radiation therapy community, which include TCP and NTCP models--is a hopeful direction \cite{nahum}. These would be particularly helpful from a treatment planning efficiency standpoint if they were embedded with uncertainties associated with organ and tumor delineation and the fact that the patient's organ/tumor geometries will change over the course of treatment. If so, such functions could be used to highlight the fact that planners often spend days on highly tailoring a treatment plan to the geometries shown in the planning CT, when the real geometries on the treatment days might be substantially different. Tackling the problem of changing geometries more head on is the subject of \emph{adaptive radiation therapy}, which is outside of the scope of this review. However, it bears mentioning that how MCO can be used in this context is an unexplored area of research.

Even without the issues of geometric uncertainties, without agreed upon objective functions the process of MCO treatment planning, as well as traditional planning, is made much less efficient. As a specific illustration, planners often find a hot spot (doses exceeding some acceptable level) in the course of working on a plan, and can spend hours trying to cool off that hot spot while not ruining other aspects of the plan. This is done on an \emph{ad hoc} basis with no way to assess the value of that hot spot (most likely, that hot spot arises in order to get complete target coverage).
With the current Pareto systems, which have at most one or two objective functions defined per contoured structure, one cannot control and assess the tradeoff of every foreseeable aspect of the dose distribution. One Pareto surface solution to this problem would be a massive Pareto surface where the dose to every voxel is considered as an independent objective. This is a new avenue for research and two questions need to be addressed. The first is, if every voxel is an objective function, does it make sense to attempt to pre-compute a Pareto surface (which would have dimension on the order of $10^6$) or is real time Pareto surface point computation based on user requests more feasible? If pre-computation is deemed a likely candidate, new theory and methods which deal with extremely large dimensional Pareto surfaces, which would naturally have many highly correlated objectives (neighboring voxels), will be needed. The result would be a Pareto based 3D dose-sculpting system. Currently, treatment planners prefer this voxel level control of the dose distributions, and they struggle to get such control from software not designed to offer it. It remains to be seen if the trend will be towards better voxel level dose sculpting tools or more reliable functions that map complex dose distributions to patient outcomes.


Another criteria that is a part of every complex treatment plan is the delivery complexity, which has a strong influence on dosimetric plan quality. Treatment planners have an incentive to keep the number of beams to a minimum in IMRT treatment planning. But, with the advent of VMAT delivery, which delivers radiation from every angle around a patient, one wonders: if one is willing to deliver a VMAT plan, why not an IMRT with 20 beams?  Theoretically, the additional treatment time needed for a 20 beam IMRT plan compared to a VMAT plan that delivers a similar dose is approximately the time it takes the gantry to rotate once fully around the patient. This is because the IMRT plan could be designed time efficiently to deliver the dose from each angle using continuous delivery, which a VMAT plan implicitly uses, thus making the only ``inefficiency'' of the IMRT plan the unused time when the gantry is moving from one delivery angle to the next. This adds up to the time for the gantry to rotate once around the patient, plus any small time that the control system needs to start and stop the gantry and to verify its position. Yet, IMRT plans with 20 beams are almost never used (see \cite{dassimrt} for an exception). There are a few reasons for this. One is that IMRT plans are often delivered beam by beam, with the radiation therapist verifying each beam before delivery. Automatic field sequencing, where all the fields are automatically delivered by the control system without therapist intervention, is an available solution to this problem (this is also called consolidated field delivery). Another issue which discourages many-beam IMRT is that quality assurance of individual treatment plans is often done on a beam by beam basis. This too is easily resolvable with existing technologies, such as QA systems that recreate the entire dose distribution by measuring the fluence from each field and summing up the results. This approach will likely become the standard for VMAT QA. Historically, many-beam IMRT plans would have also been hard to create because the additional computer memory needed to store dose matrices from all of the used beam angles. This is not a difficulty with modern computers, again evidenced by computations needed for VMAT planning.

Collectively, the historical difficulties of using many beams has led to a strong culture in today's radiotherapy planning community to use as few gantry angles (segments, etc.) as possible to create a treatment plan, and this leads to a significant amount of time grappling with the dose distribution, because often, especially for complex cases, not enough beams are used to come close an ideal Pareto optimal dose distribution. Furthermore, current planning software packages leave the treatment planners not knowing if they could do better by increasing the number of beams or changing the beam angles, etc. Using the method proposed at the end of the previous section, this difficulty is avoided.

Another impediment to the effective use of Pareto surface navigation in a clinical setting is the required change in workflow. In current standard practice in the United States, physicians are not actively involved in the planning process, they merely approve a plan as the last step. Pareto surface navigation is designed to put the tradeoff decision into the hands of the decision maker most qualified to make this decision. In the case of radiation therapy, this is the radiation oncologist. This change in workflow--having the physician partake in the planning process by navigating to a preferable tradeoff point--may prove difficult to adopt for clinics which have been doing treatment planning with the traditional workflow for decades.  If this workflow change is not implemented, a Pareto surface based system will not reach its full potential because treatment planners may navigate to a plan that fulfills the prescribed set of goals, whereas if allowed to navigate a physician might choose a plan on a different part of the Pareto surface \cite{Craft2011}.

Adopting the paradigm of physician navigation is a significant challenge for Pareto surface based treatment planning. A host of technical challenges also exist. Perhaps the largest technical challenge is achieving what is called \emph{deliverable navigation}. This refers to Pareto surface navigation such that as the planner is navigating the surface, the plan being viewed is deliverable (not for example a fluence map based IMRT or VMAT plan that still needs to be segmented, and may degrade after segmentation). The solutions for this challenge are unique to the modality. Step-and-shoot IMRT, 3D conformal therapy, VMAT, etc. will all require customized solutions. Another technical challenge, applicable to clinics continuing to use very sparse delivery techniques (very few IMRT beams, for example), will be to devise intuitive MCO methods to expose how the number of beams, segments, overall delivery time, or some other measure of treatment efficiency, influences the dose quality. In a Pareto navigation setting, this could take the form of a delivery time slider.

Finally, the application of MCO to radiation techniques requiring less advanced hardware, thus more suited to for example clinics in the developing world \cite{price2012}, such as jaws-only techniques \cite{spirou94,webbjaws,earljaws}, including perhaps jaws-only VMAT, will require creative new approaches.


\bigskip
\bibliographystyle{unsrt}
\bibliography{all}

\begin{thebibliography}{10}

\bibitem{Webb2003}
S.~Webb.
\newblock The physical basis of {IMRT} and inverse planning.
\newblock {\em British Journal of Radiology}, 76:678--689, 2003.

\bibitem{localMinima}
J.~Llacer, J.~Deasy, T.~Bortfeld, T.~Solberg, and C.~Promberger.
\newblock Absence of multiple local minima effects in intensity modulated
  optimization with dose-volume constraints.
\newblock {\em Physics in Medicine and Biology}, 48(2):183--210, 2003.

\bibitem{sens2}
M.~Alber, M.~Birkner, and F.~N{\"u}sslin.
\newblock Tools for the analysis of dose optimization: {II}. sensitivity
  analysis.
\newblock {\em Physics in Medicine and Biology}, 47(19):N265, 2002.

\bibitem{sens3}
B.~Sobotta, M.~S{\"o}hn, M.~P{\"u}tz, and M.~Alber.
\newblock Tools for the analysis of dose optimization: {III}. pointwise
  sensitivity and perturbation analysis.
\newblock {\em Physics in Medicine and Biology}, 53(22):6337, 2008.

\bibitem{nahum}
J.~Uzan and A.~Nahum.
\newblock Radiobiologically guided optimisation of the prescription dose and
  fractionation scheme in radiotherapy using {BioSuite}.
\newblock {\em British Journal of Radiology}, 85(1017):1279--1286, 2012.

\bibitem{lungheart}
E.~Huang, A.~Hope, P.~Lindsay, M.~Trovo, I.~El~Naqa, J.~Deasy, and J.~Bradley.
\newblock Heart irradiation as a risk factor for radiation pneumonitis.
\newblock {\em Acta Oncologica}, 50(1):51--60, 2011.

\bibitem{wilkens}
J.~Wilkens, J.~Alaly, K.~Zakarian, W.~Thorstad, and J.~Deasy.
\newblock {IMRT} treatment planning based on prioritizing prescription goals.
\newblock {\em Physics in Medicine and Biology}, 52(6):1675--1692, 2007.

\bibitem{falkinger}
M.~Falkinger, S.~Schell, J.~M{\"u}ller, and J.~Wilkens.
\newblock Prioritized optimization in intensity modulated proton therapy.
\newblock {\em Zeitschrift f{\"u}r Medizinische Physik}, 2011.

\bibitem{jee}
K.W. Jee, D.~McShan, and B.~Fraass.
\newblock Lexicographic ordering: intuitive multicriteria optimization for
  {IMRT}.
\newblock {\em Physics in Medicine and Biology}, 52:1845, 2007.

\bibitem{clark}
V.~Clark, Y.~Chen, J.~Wilkens, J.~Alaly, K.~Zakaryan, and J.~Deasy.
\newblock Imrt treatment planning for prostate cancer using prioritized
  prescription optimization and mean-tail-dose functions.
\newblock {\em Linear algebra and its applications}, 428(5):1345--1364, 2008.

\bibitem{craft-hmp}
D.~Craft and T.~Bortfeld.
\newblock How many plans are needed in an {IMRT} multi-objective plan database?
\newblock {\em Physics in Medicine and Biology}, 53(11):2785--2796, 2008.

\bibitem{spalke}
T.~Spalke, D.~Craft, and T.~Bortfeld.
\newblock Analyzing the main trade-offs in multiobjective radiation therapy
  treatment planning databases.
\newblock {\em Physics in Medicine and Biology}, 54(12):3741--3754, 2009.

\bibitem{ruzika}
S.~Ruzika and M.M. Wiecek.
\newblock A survey of approximation methods in multiobjective programming.
\newblock Technical report, Department of Mathematics, Univeristy of
  Kaiserslautern, 2003.

\bibitem{lo}
K.~Miettinen.
\newblock {\em Nonlinear multiobjective optimization}.
\newblock Kluwer, 2004.

\bibitem{messac}
A.~Messac, A.~Ismail-Yahaya, and C.~A. Mattson.
\newblock The normalized normal constraint method for generating the {P}areto
  frontier.
\newblock {\em Journal of the International Society of Structural and
  Multidisciplinary Optimization (ISSMO)}, 25:86--98, 2003.

\bibitem{wei-proj}
W.~Chen, D.~Craft, T.~Madden, K.~Zhang, H.~Kooy, and G.~Herman.
\newblock A fast optimization algorithm for multi-criteria intensity modulated
  proton therapy planning.
\newblock {\em Medical Physics}, 37(9):4938--4935, 2010.

\bibitem{rennen}
G.~Rennen, E.~Van Dam, and D.~Den Hertog.
\newblock Enhancement of sandwich algorithms for approximating higher
  dimensional convex {P}areto sets.
\newblock {\em {INFORMS} Journal on Computing}, 23:493--517, 2011.

\bibitem{bokrantz}
R.~Bokrantz and A.~Forsgren.
\newblock An algorithm for approximating convex pareto surfaces based on dual
  techniques.
\newblock {\em Informs Journal of Computing}, 2012.

\bibitem{monz}
M.~Monz, K-H. K\"{u}fer, T.~Bortfeld, and C.~Thieke.
\newblock Pareto navigation - algorithmic foundation of interactive
  multi-criteria {IMRT} planning.
\newblock {\em Physics in Medicine and Biology}, 53(4):985--998, 2008.

\bibitem{astroid}
H.~Kooy, B.~Clasie, H.M. Lu, T.~Madden, H.~Bentefour, N.~Depauw, J.~Adams,
  A.~Trofimov, D.~Demaret, T.~Delaney, and J.~Flanz.
\newblock A case study in proton pencil-beam scanning delivery.
\newblock {\em Int. J. Radiation Oncology Biol. Phys.}, 76(2):624--630, 2010.

\bibitem{craftmonz}
D.~Craft and M.~Monz.
\newblock Simultaneous navigation of multiple pareto surfaces, with an
  application to multicriteria {IMRT} planning with multiple beam angle
  configurations.
\newblock {\em Medical physics}, 37:736, 2010.

\bibitem{delivnav}
D.~Craft and C.~Richter.
\newblock Deliverable navigation for multicriteria step and shoot {IMRT}
  treatment planning.
\newblock {\em Physics in Medicine and Biology}, 58(1):87, 2013.

\bibitem{craft-bao}
D.~Craft.
\newblock Local beam angle optimization with linear programming and gradient
  search.
\newblock {\em Physics in Medicine and Biology}, 52(7):N127--N135, 2007.

\bibitem{evasMIPBAO}
E.~Lee, T.~Fox, and I.~Crocker.
\newblock Simultaneous beam geometry and intensity map optimization in
  intensity-modulated radiation therapy.
\newblock {\em Int. J. Radiation Oncology Biol. Phys.}, 64(1):301--320, 2006.

\bibitem{aleman2009}
D.~Aleman, H.E. Romeijn, and J.~Dempsey.
\newblock A response surface approach to beam orientation optimization in
  intensity-modulated radiation therapy treatment planning.
\newblock {\em INFORMS Journal on Computing}, 21(1):62--76, 2009.

\bibitem{BAOmdandersonWANG}
X.~Wang, X.~Zhang, L.~Dong, H.~Liu, Q.~Wu, and Mohan R.
\newblock Development of methods for beam angle optimization for {IMRT} using
  an accelerated exhaustive search strategy.
\newblock {\em Int. J. Radiation Oncology Biol. Phys.}, 60(4):1325--37, 2004.

\bibitem{icycle}
S.~Breedveld, P.~Storchi, P.~Voet, and B.~Heijmen.
\newblock {iCycle}: integrated, multicriterial beam angle, and profile
  optimization for generation of coplanar and noncoplanar {IMRT} plans.
\newblock {\em Medical Physics}, 39:951, 2012.

\bibitem{valentina}
D.~Bertsimas, V.~Cacchiani, D.~Craft, and O.~Nohadani.
\newblock A hybrid approach to beam angle optimization in intensity-modulated
  radiation therapy.
\newblock {\em Computers \& Operations Research}, 2012.

\bibitem{dao}
D.~Shepard, M.~Earl, X.~Li, S.~Naqvi, and C.~Yu.
\newblock Direct aperture optimization: A turnkey solution for step-and-shoot
  {IMRT}.
\newblock {\em Medical Physics}, 29:1007--1018, 2002.

\bibitem{dao2}
C.~Men, H.~Romeijn, Z.~Ta{\c{s}}k{\i}n, and J.~Dempsey.
\newblock An exact approach to direct aperture optimization in {IMRT} treatment
  planning.
\newblock {\em Physics in Medicine and Biology}, 52:7333, 2007.

\bibitem{daomco}
E.~Salari and J~Unkelbach.
\newblock A column-generation-based method for multi-criteria direct aperture
  optimization.
\newblock {\em Physics in Medicine and Biology}, 58(3):621, 2013.

\bibitem{vmerge}
D.~Craft, D.~McQuaid, J.~Wala, W.~Chen, E.~Salari, and T.~Bortfeld.
\newblock Multicriteria {VMAT} optimization.
\newblock {\em Medical Physics}, 39:686, 2012.

\bibitem{networkvmerge}
E.~Salari, J.~Wala, and D.~Craft.
\newblock Exploring trade-offs between {VMAT} dose quality and delivery
  efficiency using a network optimization approach.
\newblock {\em Physics in Medicine and Biology}, 57(17):5587, 2012.

\bibitem{rasmus-mco-vmat}
R.~Bokrantz.
\newblock Multicriteria optimization for volumetric-modulated arc therapy by
  decomposition into a fluence-based relaxation and a segment weight-based
  restriction.
\newblock {\em Medical Physics}, 39:6712, 2012.

\bibitem{mangasarian}
G.~Fung and O.~Mangasarian.
\newblock Equivalence of minimal l0-and lp-norm solutions of linear equalities,
  inequalities and linear programs for sufficiently small p.
\newblock {\em Journal of optimization theory and applications}, 151(1):1--10,
  2011.

\bibitem{compress}
M.~Hyder and K.~Mahata.
\newblock An approximate l0 norm minimization algorithm for compressed sensing.
\newblock In {\em Acoustics, Speech and Signal Processing, 2009. ICASSP 2009.
  IEEE International Conference on}, pages 3365--3368. IEEE, 2009.

\bibitem{zhu2009}
L.~Zhu and L.~Xing.
\newblock Search for imrt inverse plans with piecewise constant fluence maps
  using compressed sensing techniques.
\newblock {\em Medical Physics}, 36(5):1895, 2009.

\bibitem{dassimrt}
H.~Kim, R.~Li, R.~Lee, T.~Goldstein, S.~Boyd, E.~Candes, and L.~Xing.
\newblock Dose optimization with first-order total-variation minimization for
  dense angularly sampled and sparse intensity modulated radiation therapy
  ({DASSIM-RT}).
\newblock {\em Medical Physics}, 39(7):4316, 2012.

\bibitem{craft-spg}
D.~Craft, P.~S\"uss, and T.~Bortfeld.
\newblock The tradeoff between treatment plan quality and required number of
  monitor units in {IMRT}.
\newblock {\em Int. J. Radiation Oncology Biol. Phys.}, 67(5):1596--1605, 2007.

\bibitem{candes}
E.~Candes, M.~Wakin, and S.~Boyd.
\newblock Enhancing sparsity by reweighted l1 minimization.
\newblock {\em Journal of Fourier Analysis and Applications}, 14(5):877--905,
  2008.

\bibitem{pugachev02}
A.~Pugachev and L.~Xing.
\newblock Incorporating prior knowledge into beam orientaton optimization in
  imrt.
\newblock {\em Int. J. Radiation Oncology Biol. Phys.}, 54(5):1565--1574, 2002.

\bibitem{Craft2011}
D.~Craft, T.~Hong, H.~Shih, and T.~Bortfeld.
\newblock Improved planning time and plan quality through multicriteria
  optimization for intensity-modulated radiotherapy.
\newblock {\em Int. J. Radiation Oncology Biol. Phys.}, 82(1):e83--90, 2012.

\bibitem{price2012}
A.~Price, P.~Ndom, E.~Atenguena, M.~Nouemssi, J.~Pierre, and R.~Ryder.
\newblock Cancer care challenges in developing countries.
\newblock {\em Cancer}, 118(14):3627--3635, 2012.

\bibitem{spirou94}
S.~Spirou and C.~Chui.
\newblock Generation of arbitrary intensity profiles by dynamic jaws or
  multileaf collimators.
\newblock {\em Medical Physics}, 21(7):1031--1041, 1994.

\bibitem{webbjaws}
S.~Webb.
\newblock Intensity-modulated radiation therapy using only jaws and a mask.
\newblock {\em Physics in Medicine and Biology}, 47(2):257, 2002.

\bibitem{earljaws}
M.~Earl, M.~Afghan, C.~Yu, Z.~Jiang, and D.~Shepard.
\newblock Jaws-only {IMRT} using direct aperture optimization.
\newblock {\em Medical Physics}, 34:307, 2007.

\end{thebibliography}

\end{document}